\documentclass{amsart}

\usepackage{amsfonts,amssymb,verbatim,amsmath,amsthm,latexsym,textcomp,amscd}
\usepackage{latexsym,amsfonts,amssymb,epsfig,verbatim}
\usepackage{amsmath,amsthm,amssymb,latexsym,graphics,textcomp}
\usepackage{graphicx}
\usepackage{color}
\usepackage{url}

\input xy
\xyoption{all}

\begin{document}

\newtheorem{theorem}{Theorem}[section]
\newtheorem{prop}[theorem]{Proposition}
\newtheorem{lemma}[theorem]{Lemma}
\newtheorem{cor}[theorem]{Corollary}
\newtheorem{definition}[theorem]{Definition}
\newtheorem{conj}[theorem]{Conjecture}
\newtheorem{rmk}[theorem]{Remark}
\newtheorem{claim}[theorem]{Claim}
\newtheorem{defth}[theorem]{Definition-Theorem}

\newcommand{\boundary}{\partial}
\newcommand{\C}{{\mathbb C}}
\newcommand{\integers}{{\mathbb Z}}
\newcommand{\natls}{{\mathbb N}}
\newcommand{\ratls}{{\mathbb Q}}
\newcommand{\reals}{{\mathbb R}}
\newcommand{\proj}{{\mathbb P}}
\newcommand{\lhp}{{\mathbb L}}
\newcommand{\tube}{{\mathbb T}}
\newcommand{\cusp}{{\mathbb P}}
\newcommand\AAA{{\mathcal A}}
\newcommand\BB{{\mathcal B}}
\newcommand\CC{{\mathcal C}}
\newcommand\DD{{\mathcal D}}
\newcommand\EE{{\mathcal E}}
\newcommand\FF{{\mathcal F}}
\newcommand\GG{{\mathcal G}}
\newcommand\HH{{\mathcal H}}
\newcommand\II{{\mathcal I}}
\newcommand\JJ{{\mathcal J}}
\newcommand\KK{{\mathcal K}}
\newcommand\LL{{\mathcal L}}
\newcommand\MM{{\mathcal M}}
\newcommand\NN{{\mathcal N}}
\newcommand\OO{{\mathcal O}}
\newcommand\PP{{\mathcal P}}
\newcommand\QQ{{\mathcal Q}}
\newcommand\RR{{\mathcal R}}
\newcommand\SSS{{\mathcal S}}
\newcommand\TT{{\mathcal T}}
\newcommand\UU{{\mathcal U}}
\newcommand\VV{{\mathcal V}}
\newcommand\WW{{\mathcal W}}
\newcommand\XX{{\mathcal X}}
\newcommand\YY{{\mathcal Y}}
\newcommand\ZZ{{\mathcal Z}}
\newcommand\CH{{\CC\HH}}
\newcommand\PEY{{\PP\EE\YY}}
\newcommand\MF{{\MM\FF}}
\newcommand\RCT{{{\mathcal R}_{CT}}}
\newcommand\PMF{{\PP\kern-2pt\MM\FF}}
\newcommand\FL{{\FF\LL}}
\newcommand\PML{{\PP\kern-2pt\MM\LL}}
\newcommand\GL{{\GG\LL}}
\newcommand\Pol{{\mathcal P}}
\newcommand\half{{\textstyle{\frac12}}}
\newcommand\Half{{\frac12}}
\newcommand\Mod{\operatorname{Mod}}
\newcommand\Area{\operatorname{Area}}
\newcommand\ep{\epsilon}
\newcommand\hhat{\widehat}
\newcommand\Proj{{\mathbf P}}
\newcommand\U{{\mathbf U}}
 \newcommand\Hyp{{\mathbf H}}
\newcommand\D{{\mathbf D}}
\newcommand\Z{{\mathbb Z}}
\newcommand\R{{\mathbb R}}
\newcommand\Q{{\mathbb Q}}
\newcommand\E{{\mathbb E}}
\newcommand\til{\widetilde}
\newcommand\length{\operatorname{length}}
\newcommand\tr{\operatorname{tr}}
\newcommand\gesim{\succ}
\newcommand\lesim{\prec}
\newcommand\simle{\lesim}
\newcommand\simge{\gesim}
\newcommand{\simmult}{\asymp}
\newcommand{\simadd}{\mathrel{\overset{\text{\tiny $+$}}{\sim}}}
\newcommand{\ssm}{\setminus}
\newcommand{\diam}{\operatorname{diam}}
\newcommand{\pair}[1]{\langle #1\rangle}
\newcommand{\T}{{\mathbf T}}
\newcommand{\inj}{\operatorname{inj}}
\newcommand{\pleat}{\operatorname{\mathbf{pleat}}}
\newcommand{\short}{\operatorname{\mathbf{short}}}
\newcommand{\vertices}{\operatorname{vert}}
\newcommand{\collar}{\operatorname{\mathbf{collar}}}
\newcommand{\bcollar}{\operatorname{\overline{\mathbf{collar}}}}
\newcommand{\I}{{\mathbf I}}
\newcommand{\tprec}{\prec_t}
\newcommand{\fprec}{\prec_f}
\newcommand{\bprec}{\prec_b}
\newcommand{\pprec}{\prec_p}
\newcommand{\ppreceq}{\preceq_p}
\newcommand{\sprec}{\prec_s}
\newcommand{\cpreceq}{\preceq_c}
\newcommand{\cprec}{\prec_c}
\newcommand{\topprec}{\prec_{\rm top}}
\newcommand{\Topprec}{\prec_{\rm TOP}}
\newcommand{\fsub}{\mathrel{\scriptstyle\searrow}}
\newcommand{\bsub}{\mathrel{\scriptstyle\swarrow}}
\newcommand{\fsubd}{\mathrel{{\scriptstyle\searrow}\kern-1ex^d\kern0.5ex}}
\newcommand{\bsubd}{\mathrel{{\scriptstyle\swarrow}\kern-1.6ex^d\kern0.8ex}}
\newcommand{\fsubeq}{\mathrel{\raise-.7ex\hbox{$\overset{\searrow}{=}$}}}
\newcommand{\bsubeq}{\mathrel{\raise-.7ex\hbox{$\overset{\swarrow}{=}$}}}
\newcommand{\tw}{\operatorname{tw}}
\newcommand{\base}{\operatorname{base}}
\newcommand{\trans}{\operatorname{trans}}
\newcommand{\rest}{|_}
\newcommand{\bbar}{\overline}
\newcommand{\UML}{\operatorname{\UU\MM\LL}}
\newcommand{\EL}{\mathcal{EL}}
\newcommand{\tsum}{\sideset{}{'}\sum}
\newcommand{\tsh}[1]{\left\{\kern-.9ex\left\{#1\right\}\kern-.9ex\right\}}
\newcommand{\Tsh}[2]{\tsh{#2}_{#1}}
\newcommand{\qeq}{\mathrel{\approx}}
\newcommand{\Qeq}[1]{\mathrel{\approx_{#1}}}
\newcommand{\qle}{\lesssim}
\newcommand{\Qle}[1]{\mathrel{\lesssim_{#1}}}
\newcommand{\simp}{\operatorname{simp}}
\newcommand{\vsucc}{\operatorname{succ}}
\newcommand{\vpred}{\operatorname{pred}}
\newcommand\fhalf[1]{\overrightarrow {#1}}
\newcommand\bhalf[1]{\overleftarrow {#1}}
\newcommand\sleft{_{\text{left}}}
\newcommand\sright{_{\text{right}}}
\newcommand\sbtop{_{\text{top}}}
\newcommand\sbot{_{\text{bot}}}
\newcommand\sll{_{\mathbf l}}
\newcommand\srr{_{\mathbf r}}
\newcommand\geod{\operatorname{\mathbf g}}
\newcommand\mtorus[1]{\boundary U(#1)}
\newcommand\A{\mathbf A}
\newcommand\Aleft[1]{\A\sleft(#1)}
\newcommand\Aright[1]{\A\sright(#1)}
\newcommand\Atop[1]{\A\sbtop(#1)}
\newcommand\Abot[1]{\A\sbot(#1)}
\newcommand\boundvert{{\boundary_{||}}}
\newcommand\storus[1]{U(#1)}
\newcommand\Momega{\omega_M}
\newcommand\nomega{\omega_\nu}
\newcommand\twist{\operatorname{tw}}
\newcommand\modl{M_\nu}
\newcommand\MT{{\mathbb T}}
\newcommand\Teich{{\mathcal T}}
\renewcommand{\Re}{\operatorname{Re}}
\renewcommand{\Im}{\operatorname{Im}}

\title{Uniformization of simply connected finite type Log-riemann surfaces}

\author{Kingshook Biswas}
\address{RKM Vivekananda University, Belur Math, WB-711 202, India}

\author{Ricardo Perez-Marco}
\address{CNRS, LAGA, UMR 7539, Universit\'e
Paris 13, Villetaneuse, France}

\begin{abstract}
We consider simply connected log-Riemann surfaces with a finite
number of infinite order ramification points. We prove that these
surfaces are parabolic with uniformizations given by entire
functions of the form $F(z) = \int Q(z) e^{P(z)} \ dz$ where $P,
Q$ are polynomials of degrees equal to the number of infinite and
finite order ramification points respectively.
\end{abstract}

\bigskip

\maketitle

\tableofcontents

\section{Introduction}

\medskip

In \cite{bipm1} we defined the notion of log-Riemann surface, as a
Riemann surface $\SSS$ equipped with a local diffeomorphism $\pi :
\SSS \to \C$ such that the set of points $\RR$ added in the
completion ${\SSS}^* = \SSS \sqcup \RR$ of $\SSS$ with respect to
the flat metric on $\SSS$ induced by $\pi$ is discrete. The
mapping $\pi$ extends to the points $p \in \RR$, and is a covering
of a punctured neighbourhood of $p$ onto a punctured disk in $\C$;
the point $p$ is called a ramification point of $\SSS$ of order
equal to the degree of the covering $\pi$ near $p$. The finite
order ramification points may be added to $\SSS$ to give a Riemann
surface $\SSS^\times$, called the finite completion of $\SSS$. In
this article we are interested in log-Riemann surfaces of {\it
finite type}, i.e. those with finitely many ramification points and finitely
generated fundamental group, in particular simply
connected log-Riemann surfaces of finite type. We prove the following:

\medskip

\begin{theorem} \label{uniformthm} Let $\SSS$ be a log-Riemann
surface with $d_1 < +\infty$ infinite order ramification points
and $d_2 < +\infty$ finite order ramification points (counted with
multiplicity), such that the finite completion $\SSS^\times$ is
simply connected. Then $\SSS$ is biholomorphic to $\C$ and the
uniformization $\tilde{F} : \C \to \SSS^\times$ is given by an
entire function $F = \pi \circ \tilde{F}$ of the form $F(z) = \int
Q(z) e^{P(z)} dz$ where $P, Q$ are polynomials of degrees $d_1,
d_2$ respectively.
\end{theorem}

\medskip

Conversely we have:

\begin{theorem} \label{conversethm} Let $P, Q \in \C[z]$ be
polynomials of degrees $d_1, d_2$ and $F$ an entire function of
the form $F(z) = \int Q(z) e^{P(z)} dz$. Then there exists a
log-Riemann surface $\SSS$ with $d_1$ infinite order ramification
points and $d_2$ finite order ramification points (counted with
multiplicity) such that $F$ lifts to a biholomorphism $\tilde{F} :
\C \to \SSS^\times$.
\end{theorem}

\medskip

The entire functions of the above form were
first studied by Nevanlinna \cite{nevanlinna1},
who essentially proved Theorem \ref{uniformthm}, although his proof is
in the classical language. The uniformization theorem was also
rediscovered by M. Taniguchi \cite{taniguchi1} in the form of a
representation theorem for a class of entire functions defined by
him called "structurally finite entire functions". The techniques
we use are very different and adapted to the more general context of
log-Riemann surfaces. In a forthcoming article \cite{bipmlogrs2}
we use these techniques to generalize the above theorems
to a correspondence between higher genus finite type
log-Riemann surfaces and holomorphic differentials on
punctured Riemann surfaces with isolated singularities of
"exponential type" at the punctures (locally of the form $g e^h dz$ where
$g, h$ are germs meromorphic at the puncture).

\medskip

The proof of Theorem \ref{uniformthm} proceeds in outline as
follows: we approximate $\SSS$ by simply connected log-Riemann
surfaces $\SSS^\times_n$ with finitely many ramification points of
finite orders such that $d_1$ ramification points of
$\SSS^\times_n$ converge to infinite order ramification points.
The surfaces $\SSS_n$ converge to $\SSS$ in the sense of
Caratheodory (as defined in \cite{bipm1}) and by the Caratheodory
convergence theorem proved in \cite{bipm1}, the uniformizations
$\tilde{F}_n$ of $\SSS_n$ converge to the uniformization
$\tilde{F}$ of $\SSS$. The uniformizations $\tilde{F}_n$ are the
lifts of polynomials $F_n = \pi_n \circ \tilde{F}_n$, such that
the nonlinearities $G_n = F''_n/F'_n$ are rational functions of
uniformly bounded degree with simple poles at the critical points
of $F_n$. As these critical points go to infinity as $n \to
\infty$, the nonlinearity of the function $F = \pi \circ
\tilde{F}$ is a polynomial, from which it follows that $F$ is of
the form $\int Q(z) e^{P(z)} dz$.

\medskip

To prove Theorem \ref{conversethm} we use the converse of
Caratheodory convergence theorem: we approximate $F = \int Q(z)
e^{P(z)} dz$ by polynomials $F_n = \int Q(z) (1 +
\frac{P(z)}{n})^n dz$. The polynomials $F_n$ define log-Riemann
surfaces $\SSS_n$ which then converge in the sense of Caratheodory
to a log-Riemann surface $\SSS$ defined by $F$, and a study of the
log-Riemann surfaces $\SSS_n$ shows that the log-Riemann surface
$\SSS$ has $d_1$ infinite order ramification points and $d_2$
finite order ramification points (counted with multiplicity).

\medskip

We develop the tools necessary for the proofs in the following
sections. We first describe a "cell decomposition" for log-Riemann
surfaces, which allows one to approximate finite type log-Riemann
surfaces by log-Riemann surfaces with finitely many ramification
points of finite order. The cell decomposition allows us to read
the fundamental group of a log-Riemann surface from an associated
graph, and to prove a parabolicity criterion for simply connected
log-Riemann surfaces which in particular implies that the
log-Riemann surfaces $\SSS$ and $\SSS_n$ considered in the proof
of Theorem \ref{uniformthm} are parabolic.

\bigskip

\section{Cell decompositions of log-Riemann surfaces}

\medskip

We recall that a log-Riemann surface $(\SSS, \pi)$ comes equipped
with a path metric $d$ induced by the flat metric $|d\pi|$. Any
simple arc $(\gamma(t))_{t \in I}$ in $\SSS$ which is the lift of a straight line segment
in $\C$ is a geodesic segment in $\SSS$; we call such arc {\it
unbroken geodesic segments}. Note that an unbroken geodesic
segment is maximal if and only if, as $t$ tends to an endpoint of
$I$ not in $I$, either $\gamma(t)$ tends to infinity, or $\gamma(t) \to p \in
\RR$.

\bigskip

\subsection{Decomposition into stars}

\medskip

Let $w_0 \in \SSS$. Given an angle $\theta \in
\reals/2\pi\integers$, for some $0 < \rho(w_0,\theta) \leq
+\infty$, there is a unique maximal unbroken geodesic
segment $\gamma(w_0, \theta) : [0, \rho(w_0, \theta)) \to \SSS$ starting
at $w_0$ which is the lift of the line segment $\{ \pi(w_0)
+ t e^{i\theta} : 0 \leq t < \rho(w_0,\theta) \}$, such that
$\gamma(w_0, \theta)(t) \to w^* \in \RR$ if $\rho(w_0, \theta) <
+\infty$.

\medskip

\begin{definition} The star of $w_0 \in \SSS$ is the union of all
maximal unbroken geodesics starting at $w_0$, $$V(w_0) := \bigcup_{\theta \in
\reals/2\pi\integers} \gamma(w_0, \theta)$$.

Similarly we also define for a ramification point $w^*$ of order $n \leq +\infty$
the star $V(w^*)$ as the union of all maximal unbroken
geodesics $\gamma(w^*, \theta)$ starting from $w^*$, where the angle $\theta \in [-n\pi,
n\pi)$:
$$
V(w^*) := \{ \gamma(w^*, \theta)(t) : 0 \leq t < \rho(w^*, \theta) , -n\pi \leq \theta \leq n\pi \}
$$
\end{definition}

\medskip

\begin{prop} \label{stars} For $w_0 \in \SSS$ the star $V(w_0)$ is a simply connected open subset
of $\SSS$. The boundary $\partial V(w_0) \subset \SSS$ is a disjoint union of maximal unbroken geodesic
segments in $\SSS$.
\end{prop}

\medskip

\noindent{\bf Proof:} Since $\RR$ is closed, the function $\rho(w_0, \theta)$ is upper semi-continuous in
$\theta$, from which it follows easily that $V(w_0)$ is open.
Moreover $\pi$ is injective on each $\gamma(w_0, \theta)$, hence is a diffeomorphism from $V(w_0)$
onto its image $\C - F$, where $F$ is the disjoint union of closed line segments
$\{ \pi(w_0) + t e^{i\theta} : \rho(w_0,\theta) <
+\infty, t \geq \rho(w_0, \theta) \}$; clearly $\C - F$ is simply connected. By
continuity of $\pi$, each component $C$ of $\partial V(w_0)$ is
contained in $\pi^{-1}(\gamma)$ for some segment $\gamma$ in $F$,
hence is an unbroken geodesic segment $(\alpha(t))_{t \in I}$. Since $C$ is closed
in $\SSS$, $C$ must be maximal. $\diamond$

\medskip

The set of ramification points $\RR$ is discrete, hence countable.
Let $L \supset \pi(\RR)$ be the union in $\C$ of all straight
lines joining points of $\pi(\RR)$. Then $\C - L$ is dense in $\C$.
By a {\it generic fiber} we mean a
fiber $\pi^{-1}(z_0) = \{ w_i \}$ of $\pi$ such that $z_0 \in \C - L$.

\medskip

\begin{prop} Let $\{w_i\}$ be a generic fiber. Then:

\noindent{(1)} The stars $\{ V(w_i) \}$ are disjoint.

\noindent{(2)} The connected components of the stars $\partial V(w_i)$ are
geodesic rays $\gamma : (0, +\infty) \to \SSS$ such that
$\gamma(t) \to w^* \in \RR$ as $t \to 0$, $\gamma(t) \to \infty$
as $t \to \infty$.

\noindent{(3)} The union of the stars is dense in $\SSS$:
$$
\SSS = \overline{\bigcup_i V(w_i)} = \bigcup_i \overline{V(w_i)}
$$
\end{prop}

\medskip

\noindent{\bf Proof:} (1): If $w \in V(w_i) \cap V(w_j)$ then the geodesic segments
from $w$ to $w_i,w_j$ are lifts of $[\pi(w), z_0]$, so by uniqueness of lifts ($\pi$
is a local diffeomorphism) $w_i = w_j$.

\medskip

(2): By the previous Proposition, each component of $\partial
V(w_i)$ is a maximal unbroken geodesic segment $\gamma : (0, r) \to
\SSS$ with $\lim_{t \to 0} \gamma(t) = w^* \in \RR$ where $w^*$ is a ramification
point such that $\pi(\gamma)$ is a straight line segment contained in the straight line
through $\pi(w_i)$ and $\pi(w^*)$. If $r < +\infty$ then
$\gamma(t) \to w^*_1 \in \RR$ as $t \to r$, so $\pi(w_i)$ must lie
on the straight line through $\pi(w^*), \pi(w^*_1)$, contradicting
the fact that $\{w_i\}$ is a generic fiber. Hence $r = +\infty$.

\medskip

(3): Given $p \in \SSS$, if
$\pi(p) \neq z_0$, take a path $(p(t))_{0 < t < \epsilon} \subset \SSS$
converging to $p$ as $t \to 0$ such that the line segments
$[\pi(p(t)), z_0]$ make distinct angles at $z_0$,
then the discreteness of $\RR$ implies that for $t$ small enough
these line segments admit lifts; again by discreteness of $\RR$ for some $i$ we have
$p(t) \in V(w_i)$ for all $t$ small, and $p \in \overline{V(w_i)}$. $\diamond$

\medskip

It is easy to see that for $w_i \neq w_j$, the components of $\partial V(w_i), \partial V(w_j)$ are
either disjoint or equal, and each component can belong to at most two such stars.
The above Propositions hence give a cell
decomposition of $\SSS$ into cells $V(w_i)$ glued along boundary
arcs $\gamma \subset \partial V(w_i), \partial V(w_j)$.

\bigskip

\subsection{The skeleton and fundamental group}

\medskip

Let $\pi^{-1}(z_0) = \{w_i\}$ be a generic fiber. The $1$-skeleton of the
cell decomposition into stars gives an associated graph:

\medskip

\begin{definition} The {\it skeleton} $\Gamma(\SSS, z_0)$ is the graph with
vertices given by the stars $V(w_i)$, and an edge between $V(w_i)$
and $V(w_j)$ for each connected component $\gamma$ of
$\partial V(w_i) \cap \partial V(w_j)$. Each edge corresponds to a geodesic ray
$\gamma : (0, +\infty) \to \SSS$ starting at a ramification point.
This gives us a map from edges to ramification points,
$\hbox{foot} : \gamma \mapsto \hbox{foot}(\gamma)
:= \lim_{t \to 0} \gamma(t) \in \RR \in \overline{V(w_i)} \cap
\overline{V(w_j)}$.

For $w^* \in \RR$ we let $C(w^*) = \{ \gamma : \hbox{foot}(\gamma)
= w^* \}$.
\end{definition}

\medskip

We omit the proof of the following proposition which is
straightforward:

\medskip

\begin{prop} If $w^*$ is of finite order $n$ then
$C(w^*) = (\gamma_i)_{1 \leq i \leq n}$ is a cycle of edges in $\Gamma$ of
length $n$. If $w^*$ is of infinite order then $C(w^*) = (\gamma_i)_{i \in \integers}$
is a bi-infinite path of edges in $\Gamma$.
\end{prop}

\medskip

We can compute the fundamental group of a log-Riemann surface from
its skeleton:

\medskip

\begin{prop} \label{retract1} The log-Riemann surface $\SSS$ deformation retracts onto $\Gamma(\SSS,z_0)$.
In particular $\pi_1(\SSS) = \pi_1(\Gamma(\SSS,z_0))$.
\end{prop}

\medskip

\noindent{\bf Proof:} Let $\partial V(w_i) =
\sqcup_{k \in J_i} \gamma_{ik}$ be the decomposition of $\partial V(w_i)$
into its connected components. Choose points $v_{ik} \in
\gamma_{ik}$, satisfying $v_{ik} = v_{jl}$ if $\gamma_{ik} =
\gamma_{kl}$. Choose simple arcs $\alpha_{ik}, k \in J_i$, joining
$w_i$ to $v_{ik}$ within $V(w_i)$, with $\alpha_{ik} \cap \alpha_{ik'} =
\{w_i\}$. Then $\overline{V(w_i)}$ deformation retracts onto the
union of the arcs $\alpha_{ik}$; moreover for $i, j \in I$
we can choose the retractions compatibly on arcs $\gamma \subset
\partial V(w_i) \cap \partial V(w_j)$, giving a retraction of
$\SSS$ onto 
the union of all arcs
$\alpha_{ik}, i \in I, k \in J_k$, which is homeomorphic to
$\Gamma(\SSS, z_0)$. $\diamond$

\medskip

The relation of $\Gamma(\SSS, z_0)$ to the finitely completed
log-Riemann surface $\SSS^{\times}$ is as follows:

\medskip

\begin{definition} The finitely completed skeleton
$\Gamma^{\times}(\SSS, z_0)$ is the graph obtained from
$\Gamma(\SSS, z_0)$ as follows: for each finite order ramification point
$w^*$, add a vertex $v = v(w^*)$ to $\Gamma(\SSS, z_0)$, remove all edges in
the cycle $C(w^*)$ and add an edge from $v_i$ to $v$ for each vertex $v_i$ in the
cycle $C(w^*)$.
\end{definition}

Then as above we have:

\medskip

\begin{prop} \label{retract2} The finitely completed log-Riemann surface
$\SSS^{\times}$ deformation retracts onto the finitely completed
skeleton $\Gamma^{\times}(\SSS, z_0)$.
\end{prop}

\medskip

\noindent{\bf Proof:} Let $w^*$ be a finite order ramification point.
Observe that in the proof of the previous
Proposition, for $\gamma = \gamma_{ik}$ an edge in $C(w^*)$,
in the finitely completed log-Riemann surface the arc $\alpha_{ik}$
can be be homotoped to an arc $\tilde{\alpha}_{ik}$ from $w_i$ to $w^*$.
Then $\SSS^{\times}$ deformation retracts onto the union of the
arcs $\alpha_{ik}, \tilde{\alpha}_{ik}$ which is homeomorphic to
$\Gamma^{\times}(\SSS, z_0)$. $\diamond$

\medskip

Given a graph $\Gamma$ satisfying certain compatibility conditions along
with the information of the locations of the ramification points, we
can also construct an associated log-Riemann surface $\SSS$ with skeleton
$\Gamma$:

\medskip

\begin{prop} \label{logconstruct} Let $\Gamma = (V, E)$ be a connected
graph with countable vertex and edge sets and a map $\hbox{foot} : E \to \C$.
For each vertex $v$ let $E_v$ be the set of edges with a vertex at $v$ and
let $R_v = \hbox{foot}(E_v)$. Assume
that the following hold:

\smallskip

\noindent{(1)} The image $\hbox{foot}(E) \subset \C$ is discrete.

\noindent{(2)} For all vertices $v$ and points $z \in R_v$, the
intersection $\hbox{foot}^{-1}(z) \cap E_v$ has
exactly two edges, labelled $\{e_z(v,+), e_z(v,-)\}$.

\noindent{(3)} For an edge $e$ between vertices $v, v'$
with $\hbox{foot}(e) = z$, either $e = e_z(v, +)
= e_z(v', -)$ or $e = e_z(v, -)
= e_z(v', +)$.

\smallskip

Then there exists a log-Riemann surface $\SSS$ with skeleton
$\Gamma(\SSS, z_0) = \Gamma$ for some $z_0 \in \C$.

\end{prop}

\medskip

\noindent{\bf Proof:} Let $L \subset \C$ be the union of all straight lines
through pairs of points in $\hbox{foot}(E)$, and let $z_0 \in \C -
L$. For each vertex $v$ of $\Gamma$, let $L_v$ be the union of
the half-lines $l_z$ starting at points $z \in R_v$
with direction $z - z_0$. By assumption (1) this collection of half-lines is
locally finite. Let $U_v$ be the domain $\C - L_v$. Equip $U_v$
with the path metric $d(a,b) = \inf_{\beta} \int_{\beta} |dz|$
(infimum taken over all rectifiable paths $\beta$ joining $a$ and
$b$). Then the metric completion $U^*_v$ of $U_v$ is given by adjoining
for each $z \in R_v$ two copies of $l_z$ (the two 'sides' of the
slit $l_z$) intersecting at a point $z_v$, which we denote
by
$$
U^*_v = U_v \bigsqcup_{z \in R_v} (l_z(v,+) \cup l_z(v,-))
$$
where we take $l_z(v, +)$ to be the 'upper side' and $l_z(v, -)$
the 'lower' side (so $z \to l_z(v,+)$ if $z \to l_z$ in $U_v$ with $\arg(z - z_0)$
increasing and $z \to l_z(v,-)$ if $z \to l_z$ in $U_v$ with $\arg(z - z_0)$
decreasing). The inclusion of $U_v$ in $\C$ extends to a local
isometry $\pi_v : U^*_v \to \C$ with $\pi_v(l_z(v,+)) = \pi_v(l_z(v, -)) = l_z$.

\medskip

Let $\SSS^*$ be
$$
\SSS^* = \bigsqcup_{v \in V} U^*_v / \sim
$$
with the following identifications: for each edge $e$ with vertices $v, v'$
and $\hbox{foot}(\gamma) = z$, if $e = e_z(v, +) = e_z(v', -)$ we paste
isometrically the half-lines $l_z(v,+), l_z(v', -)$, otherwise we
paste isometrically $l_z(v, -), l_z(v',+)$.
The identifications are compatible with the maps $\pi_v$, giving a
a map $\pi : \SSS^* \to \C$. We let $\RR \subset \SSS^*$ be the
subset corresponding to the points $\{ z_v \}$ and $\SSS = \SSS^* - \RR$.

\medskip

Since $\pi(\RR) = \hbox{foot}(E)$ is discrete, the set $\RR$ is discrete. Moreover
$\pi$ restricted to $\SSS$ is a local isometry, and the completion
of $\SSS$ with respect to the induced path metric is precisely
$\SSS^*$, hence $\SSS$ is a log-Riemann surface. The fiber
$\pi^{-1}(z_0)$ is generic since $z_0 \in \C - L$. The stars with
respect to this fiber are precisely the open subsets $U_v \subset
\SSS$. For any star $U_v$ its closure in $\SSS^*$ is the image of $U^*_v$ in
$\SSS^*$. For vertices $v, v'$, according to the above identifications between
$U^*_{v}, U^*_{v'}$ in $\SSS^*$, each component of
$\partial U_{v} \cap \partial U_{v'}$ (if non-empty) is a half-line $l$ arising from an edge
$e$ between $v_1,v_2$, of either the form $l = l_z(v, +) =
l_z(v', -)$ or $l = l_z(v,-) = l_z(v',+)$. It follows that
$\Gamma(\SSS, z_0) = \Gamma$. $\diamond$

\bigskip

\subsection{Truncation and approximation by finite sheeted
surfaces}

\medskip

We can use the decomposition into stars to approximate any
log-Riemann surface by finite sheeted log-Riemann surfaces by "truncating"
infinite order ramification points to finite order ramification
points. More precisely we have:
%
%

\medskip

\begin{theorem} \label{truncation} Let $(\SSS, p)$ be a pointed log-Riemann
surface. Then:

\smallskip

\noindent (1) There exists a sequence of pointed log-Riemann surfaces
$(\SSS_n, p_n)$ converging to $(\SSS, p)$ in the Caratheodory
topology such that each $\SSS_n$ has only finitely many ramification
points all of finite order.

\smallskip

\noindent (2) If $\SSS^{\times}$ is simply connected then all the
surfaces $\SSS^{\times}_n$ are simply connected.
\end{theorem}

\medskip

We recall the definition of convergence of log-Riemann surfaces in
the Caratheodory topology from \cite{bipm1}: $(\SSS_n, p_n) \to
(\SSS, p)$ if for any compact $K \subset \SSS$ containing $p$
there exists $N = N(K) \geq 1$ such that for all $n \geq N$ there
is an isometric embedding $\iota_{n,K}$ of $K$ into $\SSS_n$
mapping $p$ to $p_n$ which is a translation in the charts $\pi,
\pi_n$ on $\SSS, \SSS_n$.

\medskip

\noindent{\bf Proof of Theorem \ref{truncation}:} (1): Since the generic fibers
are dense in $\SSS$ we may assume without loss of generality that $p = w_0$ lies in
a generic fiber $\{w_i\} = \pi^{-1}(z_0)$. Let $V_i = V(w_i)$ be the corresponding
stars and $\Gamma = \Gamma(\SSS, z_0)$ the associated skeleton, equipped
with the graph metric $d_{\Gamma}$ (where each edge has
length $1$). For any star $V_i$ and $R > 0$, the set $\overline{V_i} \cap \overline{B(w_i,R)}$
is compact, so it contains at most finitely many ramification
points. It follows that the collection of edges
$$
\mathcal{E}(V_i, R) := \{ \gamma : \gamma \hbox{ is an edge with a
vertex at } V_i, \hbox{foot}(\gamma) \in \overline{B(w_i,R)} \}
$$
is finite, and hence so is the corresponding collection of
vertices
$$
\mathcal{V}(V_i, R) := \{ V_j : \gamma \in \mathcal{E}(V_i,R)
\hbox{ is an edge between } V_i, V_j \}.
$$
For $n \geq 1$ we define collections of edges and vertices
$(\mathcal{E}_{n,k})_{1 \leq k \leq n}, (\mathcal{V}_{n,k})_{1 \leq k \leq n}$ as follows:

\medskip

We let $\mathcal{E}_{n,1} = \mathcal{E}(V_0, n), \mathcal{V}_{n,1} = \mathcal{V}(V_0, n)$ and for $1 < k
\leq n$,
\begin{align*}
\mathcal{E}_{n,k} & := \bigcup_{V_i \in \mathcal{V}_{n,k-1}}
\mathcal{E}(V_i, n) \\
\mathcal{V}_{n,k} & := \bigcup_{V_i \in \mathcal{V}_{n,k-1}}
\mathcal{V}(V_i,n) \\
\end{align*}

\medskip

This gives us finite connected subgraphs $\Gamma_n = (\mathcal{V}_{n,n},
\mathcal{E}_{n,n})$ of $\Gamma$ increasing to $\Gamma$. Let
$$
\hat{\SSS}_n = \bigcup_{V \in \mathcal{V}_{n,n}} \overline{V} \subset
\SSS^*
$$
be the corresponding union of stars in $\SSS^*$. It is a Riemann surface with
boundary, each boundary component being an edge $\gamma$ of $\Gamma_n$. We
paste appropriate boundary components isometrically to obtain a Riemann surface without
boundary ${\SSS}_n = \hat{{\SSS}_n} / \thicksim$ as follows:

\medskip

We let $\RR_n$ be the set of ramification points $\{
\hbox{foot}(\gamma) : \gamma \in \mathcal{E}_{n,n} \}$. For $w^* \in \RR_n$ we let
$\Gamma_n(w^*)$ be the subgraph of $\Gamma_n$ consisting of vertices $V_i$ and edges $\gamma$
such that $w^* = \hbox{foot}(\gamma) \in \overline{V_i}$. Two cases arise:

\medskip

\noindent{(i)} The ramification point $w^*$ is of finite order:
Then there are finitely many stars $V_i$ such that $w^* \in
\overline{V_i}$. If $\Gamma_n(w^*)$ does not contain all of them,
then the union of stars $\overline{V_i}, V_i \in \Gamma_n(w^*)$
has two boundary components, both of which are lifts of a
half-line in $\C$ starting at $\pi(w^*)$; in this case we can
paste the two components by an isometry which is the identity in
charts.

\medskip

\noindent{(ii)} The ramification point $w^*$ is of infinite order:
Then the union of stars $\overline{V_i}, V_i \in \Gamma_n(w^*)$
always has two boundary components, both of which are lifts of a
half-line in $\C$ starting at $\pi(w^*)$; we
paste the two components by an isometry which is the identity in
charts.

\medskip

Let $q_n : \hat{\SSS}_n \twoheadrightarrow \hat{\SSS}_n / \thicksim$ denote the
quotient of $\hat{\SSS}_n$ under the identifications made in (i), (ii). The subset
$\SSS_n := (\hat{\SSS}_n / \thicksim) - q_n(\RR_n)$ is a Riemann surface without boundary.
Since the identifications are compatible with the map $\pi$, $\pi$ induces a
map $\pi_n : \SSS_n \to \C$ which is a local diffeomorphism. The completion of $\SSS_n$ with
respect to the flat metric induced by $\pi_n$ is isometric to $\hat{\SSS}_n / \thicksim$,
so that $\SSS_n$ is a log-Riemann surface with finite ramification set
$q_n(\RR_n)$; it is clear from the construction in (i), (ii) above
that these ramification points are all of finite order. We let
$p_n = q_n(p)$.

\medskip

Any compact $K \subset \SSS$ containing $p$ can only intersect finitely
many stars $V_i$ and hence $K \subset \hat{\SSS}_n$ for $n$ large
enough. Moreover for $n$ large $K$ does not intersect the boundary
of $\hat{\SSS}_n$ (which is contained in stars going to infinity
in $\Gamma$ as $n$ goes to infinity), hence the quotient map $q_n$
isometrically embeds $K$ in $\SSS_n$. Thus $(\SSS_n, p_n)$
converges to $(\SSS, p)$ as required.

\medskip

\noindent (2): The graph $\Gamma(\SSS_n, z_0)$ can be obtained
by adding edges to the finite graph $\Gamma_n$ between certain vertices
corresponding to edges in the sets $C(w^*), w^* \in \RR_n$, to
give cycles $C(q_n(w^*))$ in $\Gamma(\SSS_n, z_0)$.
If $\SSS^{\times}$ is simply connected then by Proposition
\ref{retract2} the graph $\Gamma^{\times}(\SSS, z_0)$ is a tree.
It follows from the construction of $\Gamma^{\times}(\SSS, z_0)$
that $\pi_1(\Gamma(\SSS, z_0))$ is generated by cycles
corresponding to finite order ramification points and hence
$\pi_1(\Gamma(\SSS_n, z_0))$ is generated by the cycles $C(q_n(w^*))$.
In constructing $\Gamma^{\times}(\SSS_n, z_0)$ from
$\Gamma(\SSS_n, z_0)$ these cycles become trivial so $\pi_1(\Gamma^{\times}(\SSS_n,z_0))$
is trivial. $\diamond$


\bigskip

\subsection{Compactness for uniformly finite type log-Riemann
surfaces}

\medskip

The family of finite type log-Riemann surfaces with a given uniform
bound on the number of ramification points is compact, in the
following sense:

\medskip

\begin{theorem}\label{compactness} Let $(\SSS_n, p_n)$ be a
sequence of pointed log-Riemann surfaces with ramification sets
$\RR_n$. If for some $M, \epsilon > 0$ we have $\#\RR_n \leq M,
d(p_n, \RR_n) > \epsilon$ for all $n$ then there is a pointed
log-Riemann surface $(\SSS, p)$ with ramification set $\RR$ such
that $\#\RR \leq M$ and $(\SSS_n, p_n)$ converges to $(\SSS,p)$
along a subsequence.
\end{theorem}

%

\medskip

\noindent{\bf Proof:} Composing $\pi_n$ with a translation if necessary
we may assume $\pi_n(p_n) = 0$ for all $n$. Since $d(p_n, \RR_n) > \epsilon$ we can
change $p_n$ slightly (within the ball $B(p_n, \epsilon)$)
to assume without loss of generality that
the fiber $\pi^{-1}_n(0)$ containing $p_n$ is generic.
Let $\Gamma_n$ be the corresponding skeleton and $v_{n,0}$ the
vertex containing $p_n$.
Passing to a subsequence we may assume the projections $\pi_n(\RR_n)$ converge (in the
Hausdorff topology) to a finite set $\{ w^*_1, \dots, w^*_N\} \cup\{\infty\} \subset \hat{\C} - B(0,
\epsilon)$ (where $N \leq M$), and for all $n$ lie in small disjoint neighbourhoods $B_1, \dots,
B_N$ and $B$ of the points of $R = \{ w^*_1, \dots, w^*_N\}$ and $\infty$ respectively.

\medskip

Let $\gamma_1, \dots, \gamma_N$ be generators for the group $G = \pi_1(\C - R)$
where each $\gamma_i$ is a simple closed curve in $\C - (B \cup_i B_i)$ starting at the origin
with winding number one around $B_i$ and zero around $B_j, j \neq i$. There
is a natural action of $G$ on the vertices of $\Gamma_n$: given a vertex $v$,
let $w$ be the point of the fiber $\pi^{-1}_n(0)$ in $v$. Then any $g \in G$ has
a unique lift $\tilde{g}$ to $\SSS_n$ starting at $w$. Let $g \cdot v$ be
the vertex of $\Gamma_n$ containing the endpoint of $\tilde{g}$.

\medskip

We define a graph $\Gamma'_n = (V_n, E_n)$ as follows: the vertex set
$V_n$ is the orbit of $v_{n,0}$ under $G$. We put an edge $e$ between
distinct vertices $v,v'$ of $\Gamma'_n$ for each generator
$\gamma \in \{ \gamma^{\pm}_i, i = 1, \dots, N \}$ such that $v' =
\gamma \cdot v$. We define $\hbox{foot}_n(e) = w^*_i$ if the edge
$e$ corresponds to either of the generators $\gamma_i, \gamma^{-1}_i$. This defines a
map $\hbox{foot}_n : E_n \to R \subset \C$.

\medskip

For $v \in V_n$ let $E_v$ be the set of edges with a vertex at $v$
and $R_v = \hbox{foot}_n(E_v) \subset R$. Since $\gamma_i \cdot v =
v$ if and only if $\gamma^{-1}_i \cdot v = v$, it follows that for $z = w^*_i \in R_v$,
the intersection $\hbox{foot}^{-1}_n(z) \cap E_v$ consists of precisely the two edges
corresponding to the generators $\gamma_i, \gamma^{-1}_i$; we label these edges as
$e_z(v,+), e_z(v,-)$.

\medskip

It is easy to see that the graphs $\Gamma'_n$ satisfy the
hypotheses of Proposition \ref{logconstruct}. Since each vertex
has valence at most $2N$, the balls $B(v_{n,0}, k)$ are finite,
so we can pass to a subsequence such that
the pointed graphs $(\Gamma'_n, v_{n,0})$ converge to a limit
pointed graph $(\Gamma = (V, E), v_0)$, in the sense that for any $k \geq
1$, for all $n$ large enough there is an isomorphism $i_n$ of the
ball $B(v_0,k)$ with $B(v_{n,0}, k)$ taking $v_0$ to $v_{n,0}$. We
may also assume that the isomorphisms $i_n$ for different $n$
are compatible with the mappings $\hbox{foot}_n$ and the
labeled edges $e_z(v, +), e_z(v, -)$, thus inducing a corresponding mapping
$\hbox{foot} : E \to R \subset \C$ and a labeling of the edges of
$\Gamma$. Then the limit graph $\Gamma$ satisfies the hypotheses
of Proposition \ref{logconstruct} and we obtain a corresponding
pointed log-Riemann surface $(\SSS, p)$ ramified over the points of $R$
such that $\Gamma(\SSS, 0) = \Gamma$,
with $p$ in a generic fiber $\pi^{-1}(0)$, and the star containing
$p$ corresponding to the vertex $v_0$ of $\Gamma$. Moreover $\SSS$ has
at most $N$ ramification points. It is easy to
see that any compact $K \subset \SSS$ containing $p$ embeds
isometrically in all the log-Riemann surfaces $\SSS_n$ via an
isometry $\iota_n$ such that $\iota_n(p) = p_n, \iota'_n(p) = 1$,
hence $(\SSS_n, p_n)$ converges to $(\SSS, p)$. $\diamond$

\medskip


\bigskip

\subsection{Decomposition into Kobayashi-Nevanlinna cells}

\medskip

Let $\SSS$ be a log-Riemann surface with $\RR \neq \emptyset$. We
define a cellular decomposition of $\SSS$ due to Kobayashi \cite{kobayashi}
and Nevanlinna (\cite{nevanlinna2}  which is useful in determining
the type (parabolic or hyperbolic) of simply connected log-Riemann
surfaces.

\medskip

\begin{definition} \label{kobnevcell} Let $w^* \in \RR$. The
Kobayashi-Nevanlinna cell of $w^*$ is defined to be the set
$$
W(w^*) := \{ w \in \SSS^* | d(w,w^*) < d(w, \RR - \{w^*\}) \}
$$
\end{definition}

\medskip

\begin{prop} The Kobayashi-Nevanlinna cells satisfy:

\noindent (1) Any $w \in W(w^*)$ lies on an unbroken geodesic
$[w^*, w] \subset W(w^*)$. In particular $W(w^*) \subset V(w^*)$ is open
and path-connected.

\noindent (2) The boundary of $W(w^*)$ is a locally finite union of geodesic
segments.

\noindent (3) $\SSS = \overline{\cup_{w^* \in \RR} W(w^*)}$
\end{prop}

\medskip

\noindent{\bf Proof:} (1): For any $w \in W(w^*), w \neq w^*$, since $\RR \neq \emptyset$
there is a maximal unbroken geodesic $\gamma(w, \theta)$ converging to a point of $\RR$ at one
end, and since $w^*$ is the point in $\RR$ closest to $w$, there
must be such a geodesic $[w, w^*]$ converging to $w^*$. Moreover
for any $w' \in [w, w^*], w^*_1 \in \RR - \{w^*\}$, we have
$$
d(w^*, w') = d(w^*, w) - d(w, w') < d(w^*_1, w) - d(w, w') \leq
d(w^*_1, w')
$$
hence $[w, w^*] \subset W(w^*)$.

\medskip

\noindent (2): Let $w \in \partial W(w^*)$. By discreteness of
$\RR$ there are finitely many ramification points $w^* = w^*_1, \dots,
w^*_n$ at minimal distance $r > 0$ from $w$, and $n \geq 2$. The disc
$B(w, r)$ is a euclidean disk, with the points $w^*_i$ lying on its boundary;
the angular bisectors of the sectors formed by $[w, w^*_i], [w, w^*_{i+1}]$
then are equidistant from $w^*_i, w^*_{i+1}$ and lie in $\partial W(w^*_i) \cap \partial W(w^*_{i+1})$,
while all other points in the disk lie in $W(w^*_i)$ for some $i$.
Hence a neighbourhood of $w$ in $\partial W(w^*)$ is given either by a
geodesic segment passing through $w$ (if $n = 2$) or by two
geodesic segments meeting at $w$ (if $n > 2$).

\medskip

\noindent (3): Any $w \in \SSS$ belongs to $\overline{W(w^*)}$ for any
ramification point $w^*$ at minimal distance from $w$. $\diamond$

\bigskip

\subsection{Kobayashi-Nevanlinna parabolicity criterion}

\medskip

We consider a log-Riemann surface $\SSS$ such that the finite
completion $\SSS^\times$ is simply connected.
We will use the following theorem of Nevanlinna
(\cite{nevanlinna2} p. 317):

\begin{theorem} \label{kobnevcriterion} Let $F \subset {\SSS}^{\times}$ be a discrete set and
$U : {\SSS}^{\times} - F \to [0,+\infty)$ be a continuous function such that:

\noindent (1) $U$ is $C^1$ except on at most a family of locally
finite piecewise smooth curves.

\noindent (2) $U$ has isolated critical points.

\noindent (3) $U \to +\infty$ as $z \to F$ or as $z \to \infty$.

For $\rho > 0$ let $\Gamma_{\rho}$ be the union of the curves
where $U = \rho$, and let
$$
L(\rho) = \int_{\Gamma_{\rho}} |\hbox{grad}_z U| |dz| .
$$
where $|\hbox{grad}_z U| |dz|$ is the conformally invariant
differential given by
$\sqrt{\left(\frac{\partial U}{\partial x}\right)^2 + \left(\frac{\partial U}{\partial
y}\right)^2}|dz|$ for a local coordinate $z = x + iy$.
If the integral
$$
\int_{0}^{\infty} \frac{d\rho}{L(\rho)}
$$
is divergent then the surface ${\SSS}^{\times}$ is parabolic.
\end{theorem}

\medskip

We now define a function $U$ on $\SSS$ as follows:

\medskip

Let $\omega$ be the continuous differential $\omega := |d\arg(w -
w^*)|$, where for each $w \in \SSS$, $w^*$ is a ramification point
such that $w \in \overline{W(w^*)}$. Fix a base point $w_0 \in \SSS$ and
define $\tau : \SSS \to [0,+\infty)$ by
$$
\tau(w) := \inf \int_{w_0}^{w} \omega
$$
where the infimum is taken over all paths from $w_0$ to $w$.
We define another nonnegative continuous function $\sigma : \SSS
\to [0,+\infty)$ by
$$
\sigma(w) := |\log|w - w^*||
$$
where as before for each $w \in \SSS$ the point $w^*$ is a
ramification point such that $w \in \overline{W(w^*)}$.

\medskip

Then the sum $U = \tau + \sigma : \SSS \to \reals$ is a function
satisfying the conditions (1)-(3) of the above theorem. The map
$t = \sigma+i\tau$ gives a local holomorphic coordinate away from the
boundaries of the Kobayashi-Nevanlinna cells, for which we have $|\hbox{grad}_t U| |dt| =
\sqrt{2}|dt|$. On a level set $\Gamma_{\rho} = \{ U = \rho \}$ we
have $0 \leq \tau \leq \rho, t = (\rho - \tau) + i\tau$, so
$|\hbox{grad}_t U| |dt| = \sqrt{2}|dt| = 2|d\tau|$. For a given $\theta > 0$, the connected components
of the level set $\{ \tau(w) = \theta \}$ are Euclidean line segments which are
half-lines or intervals; let $0 \leq n(\theta) \leq \infty$ denote the number of such
line segments. Each such segment intersects $\Gamma_{\rho}$ in at
most one point; hence we obtain
$$
L(\rho) = \int_{\Gamma_{\rho}} |\hbox{grad}_t U| |dt| = 2 \int_{\Gamma_{\rho}} |d\tau|
\leq \int_{0}^{\rho} n(\theta) d\theta
$$
Using Theorem \ref{kobnevcriterion} above, we obtain the
following:

\medskip

\begin{theorem} \label{paracriterion} Let $\SSS$ be a log-Riemann
surface such that $\SSS^\times$ is simply connected. For $\theta >
0$ let $0 \leq n(\theta) \leq \infty$ denote the number of connected components of the level set
$\{ \tau(w) = \theta \}$. If the integral
$$
\int_{0}^{\infty} \frac{d\rho}{\int_{0}^{\rho} n(\theta) d\theta}
$$
is divergent then $\SSS^\times$ is biholomorphic to $\C$.
\end{theorem}

\medskip

This implies:

\begin{cor} \label{parabolic} Let $\SSS$ be a log-Riemann surface with a
finite number of ramification points such that $\SSS^{\times}$ is
simply connected. Then $\SSS$ is biholomorphic to $\C$.
\end{cor}

\medskip

\noindent{\bf Proof:} In this case the function $n(\theta)$ is
bounded above by twice the number of ramification points of
$\SSS$, so $\int_{0}^{\rho} n(\theta) d\theta \leq C\rho$ and
hence the integral in Theorem \ref{paracriterion} diverges.
$\diamond$

\bigskip

\section{Uniformization theorems}

\medskip

We can now prove Theorem \ref{uniformthm} as follows:

\medskip

\noindent{\bf Proof of Theorem \ref{uniformthm}:} Let $p \in \SSS$. Let $D_1, D_2$ be the
numbers of infinite and finite order ramification points respectively of $\SSS$. By Corollary
\ref{parabolic} the log-Riemann surface $\SSS^{\times}$ is
biholomorphic to $\C$. The approximating finitely completed log-Riemann surfaces
$\SSS^{\times}_n$ given by Theorem \ref{truncation} are also biholomorphic
to $\C$ and for $n$ large all have $D_1 + D_2$ ramification points.
Let $\tilde{F} : \C \to \SSS^{\times}$ and $\tilde{F}_n : \C \to \SSS^{\times}_n$
be corresponding normalized uniformizations such that
$\tilde{F}(0) = p, \tilde{F}'(0) = 1, \tilde{F}_n(0) = p_n, \tilde{F}'_n(0) =
1$, with inverses $G = \tilde{F}^{-1}, G_n = \tilde{F}^{-1}_n$.
By Theorem 1.2 of \cite{bipm1} the
entire functions $F_n = \pi_n \circ \tilde{F}_n$ converge
uniformly on compacts to the entire function $F = \pi \circ
\tilde{F}$. Since $\pi_n : \SSS^{\times}_n \to \C$ is finite to one, the entire
function $F_n$ has a pole at $\infty$ of order equal to the degree of $\pi_n$, and
is hence a polynomial. The nonlinearities $R_n = F''_n / F'_n$ are rational
functions whose poles are simple poles with integer residues at the
critical points of $F_n$, which are images of the
ramification points of $\SSS_n$ under $G_n$. Thus the rational functions $R_n$
are all of degree $D_1 + D_2$, converging normally to $F'' /
F'$, so $R = F''/F'$ is a rational function of degree at most
$D$.

\medskip

Each ramification point $w^*$ of $\SSS$ corresponds to a ramification point
$w^*_n$ of $\SSS_n$ of order converging to that of $w^*$.
We note that for $n$ large any compact $K \subset \SSS^{\times}$ containing
$p$ embeds into the approximating surfaces $\SSS^{\times}_n$.
Since the maps $G_n$ converge to $G$ uniformly on compacts of $\SSS^{\times}$ by
Theorem 1.1 \cite{bipm1}, the images under $G_n$ of ramification points
in $\SSS^{\times}_n$ corresponding to finite ramification points in $\SSS$ converge
to their images under $G$, giving in the limit $D_2$ simple poles of $R$, with residue at each
equal to the order of the corresponding finite ramification point of $\SSS$ minus
one.

\medskip

On the other hand the infinite order ramification points of $\SSS$
are not contained in $\SSS^{\times}$, so the images of the corresponding
ramification points in $\SSS^{\times}_n$ under $G_n$ cannot be contained
in any compact in $\C$ and hence converge to infinity. The
rational functions $R_n$ have a simple zero at infinity, and have
$D_1$ simple poles converging to infinity. Applying the Argument
Principle to a small circle around infinity it follows that $R$
has a pole of order $D_1 - 1$ at infinity.

\medskip

Thus $R$ is of the form
$$
\frac{F''}{F'} = \sum_{i = 1}^{D_2} \frac{m_i - 1}{z - z_i} +
P'(z)
$$
where $m_1, \dots, m_{D_2}$ are the orders of the finite ramification points
of $\SSS$ and $P$ is a polynomial of degree $D_1$. Integrating the above
equation gives
$$
F(z) = \pi(p) + \int_{0}^{z} (t - z_1)^{m_1 - 1} \dots (t -
z_{D_2})^{m_{D_2} - 1} e^{P(t)} dt
$$
as required. $\diamond$

\medskip

We can prove the converse using the above Theorem and the compactness Theorem.
We need a lemma:

\medskip

\begin{lemma} \label{simply connected} Let $(\SSS_n, p_n)$
converge to $(\SSS, p)$. If all the surfaces $\SSS^{\times}_n$ are
simply connected then $\SSS^{\times}$ is simply connected.
\end{lemma}

\medskip

\noindent{\bf Proof:} We may assume the points $p_n,p$ belong to generic
fibers. Let $\Gamma_n, \Gamma$ denote the corresponding skeletons.
Let $\gamma$ be a loop in $\SSS^{\times}$
based at $p$. We may homotope $\gamma$ away from the finite
ramification points to assume that $\gamma \subset \SSS$. By
Proposition \ref{retract1}, $\gamma$ corresponds to a path of
edges $\alpha = \{e_1, \dots, e_n\}$. By induction on the number of edges we
may assume that $\alpha$ is simple. If
$\hbox{foot}(\alpha) = \{w^*\}$ is a singleton then $w^*$ is
a finite ramification point and $\gamma$ is
trivial in $\SSS^{\times}$. Otherwise there are distinct
ramification points $w^*_1, w^*_2 \in \hbox{foot}(\alpha)$.
Considering the isometric embedding of $\gamma$ in
$\SSS_n$ for $n$ large gives a path $\gamma_n$ and a corresponding path of
edges $\alpha_n$; for $n$ large, it follows that
there are distinct ramification points in
$\hbox{foot}_n(\alpha_n)$, hence $\gamma_n$ is non-trivial in
$\SSS^{\times}_n$, a contradiction. $\diamond$


%

\medskip

\noindent{\bf Proof of Theorem \ref{conversethm}:} Given an entire
function $F$ with $F'(z) = Q(z) e^{P(z)}$ we can approximate it by
polynomials $F_n$ such that $F'_n(z) = Q(z)(1 + P(z)/n)^n$. Let
$Z_n  = \{ P = -n \} \cup \{ Q = 0 \} \cup  \subset \C$
be the zeroes of $F'_n$. The
pair $(\SSS_n = \C - Z_n, \pi_n = F_n : \C - Z_n \to \C)$ is a
log-Riemann surface with finite ramification set $\RR_n$ which can be naturally
identified with $Z_n$, the order of a ramification point
being the local degree of $F_n$ at the corresponding point of
$Z_n$.

\medskip

For $n$ large the surfaces $\SSS_n$ all have the same number of
ramification points $D = D_1 + D_2$ where $D_1$ is the degree of
$P$ and $D_2$ the number of distinct zeroes of $Q$. Moreover since
$F'_n$ converge uniformly on compacts, choosing a point $z_0$ such
that $Q(z_0) \neq 0$, for all $n$ large $|F'_n|$ is uniformly
bounded away from $0$ on a fixed neighbourhood of $z_0$, so
$d(z_0, \RR_n)$ is uniformly bounded away from $0$. It follows
from Theorem \ref{compactness} that $(\SSS_n, p_n = z_0)$ converge
along a subsequence to a limit log-Riemann surface $(\SSS, p)$
with finitely many ramification points such that $\pi(p) = z_0$.
Since $\SSS^{\times}_n$ is simply connected for all $n$, by the
previous Lemma $\SSS^{\times}$ is simply connected. By Theorem
\ref{parabolic}, $\SSS^{\times}$ is biholomorphic to $\C$. Let
$\tilde{F} : \C \to \SSS^{\times}$ be a normalized uniformization
such that $\tilde{F}(z_0) = p, \tilde{F}'(z_0) = F'(z_0)$. It
follows from Theorem 1.2 of \cite{bipm1} that the maps $F_n$
converge normally to $\pi \circ \tilde{F}$, so $F = \pi \circ
\tilde{F}$. Thus $F$ defines the uniformization of a simply
connected log-Riemann surface with finitely many ramification
points. The degrees of $Q, P$ relate to the numbers of finite
poles and poles at infinity respectively of the nonlinearity
$F''/F'$; the relations between the degrees of $Q,P$ and the
numbers of finite and infinite order ramification points of $\SSS$
then follow from the previous Theorem. $\diamond$

\bibliography{uniform}
\bibliographystyle{alpha}

\end{document}